\newcommand{\be}{\begin{eqnarray}}
\newcommand{\ee}{\end{eqnarray}}

\def\ni{\noindent}

\def\d{\mbox{$\delta$}}
\def\bd{\begin{displaymath}}
\def\ed{\end{displaymath}}
\def\ba#1{\begin{array}{#1}}
\def\ea{\end{array}}
\def\nn{\nonumber}
\newfont{\Bbb}{msbm10 scaled 1200}
\documentstyle[12pt]{article}
\title{
Lie Bialgebra Structures\\
on Twodimensional Galilei\\
Algebra and their Lie--Poisson\\
Counterparts.}
\author{Emil Kowalczyk\thanks{Supported by the \L\'od\'z University Grant No.487}\\
Department of Field Theory\\
University of \L\'od\'z \\
Pomorska 149/153\\
90--236 \L\'od\'z, Poland
}
\date{}
\begin{document}
\maketitle

\begin{abstract}
All bialgebra structures on twodimensional Galilei algebra
are classified. The corresponding Lie--Poisson structures on
Galilei group are found.
\end{abstract}
\newpage
\section{Introduction}
Much interest has been recently attracted to the problem of deformations
of space--time symmetries (see~[1]--[5] and references contained therein).
Most papers are devoted to the relativistic case. There is, however, a number
of papers dealing with Galilean symmetries also. These are the papers of
Firenze group~[6]--[8] as well as some others~[9]--[13].

The present paper is the first of the series devoted to systematic study
of quantum Galilei groups. We are going to study the following problems:\\
--the classification of all bialgebra structures (resp. Lie--Poisson structures)
on twodimensional Galilei algebra (resp. group);\\
--the classification of such structures in the fourdimensional case;\\
--quantization of the resulting Lie bialgebras and Poisson--Lie groups;
their structure and duality relations;\\
--nonrelativistic quantum space--times;\\
--differential calculi on quantum Galilei groups and 
corresponding space--times;\\
--representations, central extensions;\\
--physical applications;

In this first paper we classify the bialgebra structures on twodimensional
Galilei algebra. To this end we find the most general 1--cocycle $\delta$.
However, two such $\delta$'s should be viewed as equivalent if
they are related under the automorphism of the algebra. We find
the action of automorphism group in the space of 1--cocycles and 
classify the orbits. This allows to find all nonequivalent Lie--Poisson
structures on Galilei group. The whole procedure follows quite closely
the one presented in Ref.~[14] for $E(2)$ group.

As a result we find nine nonequivalent bialgebra structures--four 
families parametrized by nonnegative dimensionless parameter and
five discrete ``points''. In the general case $\delta$ is a sum
of few terms which makes necessary to introduce dimensionful 
parameters to provide the proper dimensions to all terms. They are
arbitrary but fixed and nonvanishing; varying them we obtain equivalent
bialgebras.
The number of dimensionful parameters can be  in general reduced by 
ascribing proper dimension to the Poisson brackets.

It is remarkable that only one discrete case corresponds to coboundary
\d. This is in sharp contrast with the semisimple case~[15] as well
as fourdimensional Poincare group~[16]~[17].\\
The bialgebra (or Lie--Poisson) structures appearing in the literature
fit into the scheme. For example the one from Ref.~[13] corresponds to
first case in Table~1 with $\varepsilon=0$; the bialgebra structure 
inferred from Ref.~[6]
with vanishing central element $M$ lies on the orbit described in 
fourth row of Table~1 with $\varepsilon=1$.

Actually, our results to large extent coincide with those of Ref.[18]
where the Heisenberg--Weyl group was considered; this is due to the fact
that both groups are isomorphic. However, our aim was to give a precise
classification of all structures which are nonequivalent with respect
to the action of the group of automorphisms.
\section{Twodimensional Galilei group,algebra and their
 automorphisms}
Twodimensional Galilei group is defined as the group of transformations
 of twodimensional space--time.
$$
\left.\begin{array}{rclcl}
x&\to&x^\prime&=&x+v t+a\\
t&\to&t^\prime&=&t+\tau
\end{array}\right\}\eqno (1)$$
Accordingly, the composition law reads
\setcounter{equation}{1}

\be
(\tau,v,a)\star(\tau^\prime,v^\prime,a^\prime)&=&(\tau+\tau^\prime,v+v^\prime,a+a^\prime+v\tau^\prime)
\ee

There exists the following global exponential parametrization of group
elements

\be
g&=&e^{-i\tau H}e^{iaP}e^{iv K}
\ee

\noindent which defines the generators of time translations $(V)$,
space translations $(P)$\ and boots $(K)$; their dimensions, respectively, 
are: $[H]=($time$)^{-1}, [P]=($distance$)^{-1}, [K]=($velocity$)^{-1}$.

The resulting Lie algebra reads

\be
[K,H]=iP,&[K,P]=0,&[H,P]=0
\ee

This algebra can be realized in terms of left--invariant fields to be
calculated according to the standard rules  from the composition law~(2):\\

\be
H&=&i(\frac{\partial}{\partial v}+v\frac{\partial}{\partial a})\equiv X^L_\tau\nonumber\\
P&=& -i\frac{\partial}{\partial a}\equiv X^L_a\\
K&=&-i\frac{\partial}{\partial v}\equiv X^L_v\nonumber
\ee

One can also compute the right--invariant fields

\be
X^R_\tau&=&-i\frac{\partial}{\partial\tau}\nonumber\\
X^R_a&=&i\frac{\partial}{\partial a}\\
X^R_v&=&i(\frac{\partial}{\partial v}+\tau\frac{\partial}{\partial a})\nonumber
\ee

\noindent which obey the same commutation rules.

The following observation will be important below. If we redefine \\
$H\to H/\tau_0, P\to P/v_0\tau_o, K\to K/v_0$, where $\tau_0$\ and $v_0$
are arbitrary constants of dimension time and velocity, respectively, 
we obtain the same algebra (4), the generators $H,P,K$\ being now, however, 
dimensionless.

It is very easy to describe all automorphisms of the algebra (4) 
(taken in dimensionless form).\\
Let
\setcounter{equation}{7}
$$
A=\left(\begin{array}{cc}
a&b\\
c&d
\end{array}\right)\in GL(2,\mbox{\Bbb R}),\;\; \vec{X}=\left(\begin{array}{c}
x_1\\
x_2
\end{array}\right);\eqno (\theequation)
$$
\noindent then the general automorphism $H$\ can be written in form
\setcounter{equation}{8}
$$
\left(\begin{array}{c}
K\\
H\\
P
\end{array}\right)\to
\left(\ba{c}
K^\prime\\
H^\prime\\
P^\prime
\ea\right)=\left(\ba{cc}
A&\vec{X}\\
0&detA\ea\right)\left(\ba{c}
K\\
H\\
P\ea\right) \eqno (\theequation)
$$
The group of automorphisms is therefore sixdimensional. The composition law read

\be
(A,\vec{X})\star(\vec{A}^\prime,\vec{X}^\prime)&=&(A\cdot A^\prime, \vec{AX^\prime}+(detA^\prime)\vec{X})
\ee

\section{The bialgebra structures on twodimensional Galilei algebra}
In this section we find all bialgebra structures on twodimensional Galilei
algebra. As it is well known  given any Lie algebra L the bialgebra $(L,\delta)$
is defined by a skewsymmetric cocommutator $\delta:\;L\to L\otimes L$\ such
that:\\
(i) $\delta$\ is a 1--cocycle, i.e.

\be
\delta([X,Y])&=&[\delta(X),1\otimes Y+Y\otimes 1]+[1\otimes X+X\otimes 1,\delta(Y)]\;\;X,Y\in L\nonumber
\ee

(ii) the dual map $\delta^\star:\;L^\star\otimes L^\star\to L^\star$\ defines
a Lie bracket on $L^\star$\\
Our aim is to find all bialgebra structures on the algebra (4). The general
form of $\delta$\ obeying (i) is\\

\be
\delta(P)&=&sK\wedge P-rH\wedge P\nonumber\\
\delta(H)&=&sK\wedge H+\alpha K\wedge P-\beta H\wedge P\\
\delta(K)&=&rK\wedge H+\gamma K\wedge P-\rho H\wedge P\nonumber
\ee

\noindent$s,r,\alpha,\beta,\gamma,\rho$\ being arbitrary real parameters.\\
The condition (ii) adds two further constraints\\

\be
2s\rho-r(\beta+\gamma)&=&0\\
2r\alpha-s(\beta+\gamma)&=&0\nonumber
\ee

Eqs. (10)--(11) define all bialgebra structures on twodimensional Galilei
algebra. However, it is not the end of the story. Two $\delta$'s 
which can be transformed into each other by automorphism (8) are 
considered as equivalent. Therefore, we are interested in classification of
nonequivalent bialgebra structures. To this end we find first the transformation
rules for the parameters. They read

\be
\alpha^\prime&=&\frac{1}{\Delta^2}(d^2\alpha+cd(\beta+\gamma)+c^2\rho)\nonumber\\
\beta^\prime&=&\frac{1}{\Delta^2}[(bd\alpha+ad\beta+bc\gamma+ac\rho)-(cr+ds)x_1+(ar+bs)x_2]\nonumber\\
\gamma\; ^\prime&=&\frac{1}{\Delta^2}[(bd\alpha+bc\beta+ad\gamma+ac\rho)+(cr+ds)x_1-(ar+bs)x_2]\\
\rho^\prime&=&\frac{1}{\Delta^2}(b^2\alpha+ab(\beta+\gamma)+a^2\rho)\nonumber\\
r^\prime&=&(ar+bs)\frac{1}{\Delta}\nn\\
s^\prime&=&(cr+ds)\frac{1}{\Delta},\,\; \Delta\equiv detA\equiv ad-bc\nn
\ee

We see that, apart from $\Delta$--factors, $r,s$\ transform according to the
defining representation of $GL(2,\mbox{\Bbb R}),\alpha,\beta+\gamma,\rho$\ 
form  an irreducible triplet from symmetric product of 
two defining representations
while the transformation rule for $\beta-\gamma$\ reads

\be
\beta^\prime-\gamma\; ^\prime&=&\frac{1}{\Delta}(\beta-\gamma)+\frac{2}{\Delta^2}
(-(cr+ds)x_1+(ar+bs)x_2)
\ee

In order to classify the inequivalent $\delta$'s we observe first that, as
far as $r$\ and $s$\ are concerned, there are two orbits: $r=s=0$\ or
$|r|+|s|\neq 0$. Take $r=s=0$; denoting $T_{11}=\rho, T_{22}=\alpha, 
T_{12}=\frac{1}{2}(\beta+\gamma)$\ one gets

\be
T^\prime&=&\frac{1}{(detA)^2}ATA^T\\
\beta^\prime-\gamma\; ^\prime&=&\frac{1}{\Delta}(\beta-\gamma)
\ee

Moreover, there are no further constaints on $\alpha,\beta,\gamma$\ and
$\rho$\ because ,  eqs.~(11) (which are invariant under the automorphisms~(12)~)
are obeyed automatically. 

Let us diagonalize $T$\ by choosing on appropriate orthogonal $A$. 
According to the Sylvester theorem there are the following inequivalent
possibilities for eigenvalues $\lambda_1,\lambda_2$\ of $T$:\\

\ni a)$\lambda_1>0,\lambda_2>0$\\
b)$\lambda_1<0,\lambda_2<0$\\
c)$\lambda_1>0,\lambda_2<0$\\
d)$\lambda_1>0,\lambda_2=0$\\
e)$\lambda_1<0,\lambda_2=0$\\
f)$\lambda_1=0,\lambda_2=0$\\
\vskip1ex
\noindent Taking further\\
$$
A=\left(\ba{cc}
\mu_2&0\\
0&\mu_1
\ea\right), \mu_1\mu_2\neq 0
$$
\noindent we get $\lambda_i\to\lambda_i/{\mu^2_i}$. For $\lambda_1\lambda_2\neq 0$, i.e for (a)--(c)
cases we take $\mu_i=\sqrt{|\lambda_i|}$\ enforcing $|\lambda_i|=1$.
For (a),(b) cases $T$\ is further invariant under $0(2)$\ transformations
while for the case (c)--under $0(1,1)$\ ones. Indeed, taking the determinant
of both sides of equation\\
$$
T=\frac{1}{(detA)^2}ATA^T
$$
\noindent we get
$$
(detA)^2=1,{\rm i.e.}\;\;\; detA=\pm 1$$
\noindent and
$$
T=ATA^T
$$
By choosing an appropriate sign of $detA$\ one can achieve, according to eq.~(15),\ 
$\beta^\prime-\gamma\; ^\prime \geq 0$\\

If (d) or (e) holds, $\mu_2$\ can be arbitrary while $\mu_1$\ is chosen
to enforce $|\lambda_1|=1$. Therefore, according to eq.~(15) one can take $\mu_2$\
such that $\beta^\prime-\gamma\; ^\prime=0$, respectively $\beta^\prime-\gamma\; ^\prime=2$\ 
depending on whether $\beta-\gamma=0$, respectively $\beta-\gamma\neq 0$;
the same holds for (f).

Summarizing, for $r=s=0$\ there are the following inequivalent possibilities:\\

\be
\alpha=1,& \rho=1,& \beta=-\gamma=\varepsilon \geq 0\nn\\ 
\alpha=-1,& \rho=-1,& \beta=-\gamma=\varepsilon \geq 0\nn\\ 
\alpha=-1,& \rho=1,& \beta=-\gamma=\varepsilon \geq 0\nn\\ 
\alpha=0,& \rho=1,& \beta=\gamma=0\nn\\ 
\alpha=0,& \rho=1,& \beta=-\gamma=1\nn\\ 
\alpha=0,& \rho=-1,& \beta=\gamma=0\nn\\ 
\alpha=0,& \rho=-1,& \beta=-\gamma=1\nn\\ 
\alpha=0,& \rho=0,& \beta=-\gamma=1\nn
\ee

Let us now consider the case $|r|+|s|\neq 0$. There is now one orbit for the
$(r,s)$--dublet. Therefore we can choose $r^\prime=1,\;s^\prime=0$, i.e.
\setcounter{equation}{16}
$$\left\{\ba{rcl}
ar+bs&=&\Delta \\
cr+ds&=&0 \ea\right.\eqno (\theequation)
$$
\noindent Solving eqs.(16) one obtains

\be
c=-s&&d=r
\ee

\ni This, together, with eqs.(11)~and~(12) implies

\be
\alpha^\prime=0&&\beta^\prime+\gamma'=0\nn
\ee

\ni On the other hand eq.(13) gives

\be
\beta'-\gamma'&=&\frac{1}{\Delta}((\beta-\gamma)+2x_2)
\ee

Therefore, taking $x_2=-(\frac{\beta-\gamma}{2})$\ and c,d as determined
by eq.(17) we arrive at $\alpha'=\beta'=\gamma'=s'=0, r'=1$; it remains to find
the possible values of $\rho'$. Now, eqs.(11) imply also

\be
r((\beta+\gamma)^2-4\alpha\rho)=0,&s((\beta+\gamma)^2-4\alpha\rho)=0&\; \nn
\ee

\noindent or, due to $|r|+|s|\neq 0$,

\be
(\beta+\gamma)^2&=&4\alpha\rho 
\ee

The quadratic form

\be
\alpha b^2+\rho a^2+(\beta+\gamma)ab& & \nn
\ee

\ni entering the expression for $\rho'$, is semidefinite. It vanishes provided

\be
2\rho a+(\beta+\gamma)b&=&0
\ee

\ni which has no common solution with $ra+sb=\Delta$\ $(2\rho s-r(\beta+\gamma)=0!)$
unless $\rho=(\beta+\gamma)=0$. If the latter holds we can take $\rho'=0$.
If not, $\rho'>0$, resp. $\rho'<0$\ if $\rho>0$, resp. $\rho<0$. Nothing
more can be done because if we start with $\alpha=\beta+\gamma=s=0$, $r=1$,
then, according to the eqs.(12)~and~(17), $\Delta=a$\ and, consequently,
$\rho=\rho'$.

Accordingly, for $|r|+|s|\neq 0$\ we have the following canonical position:\\
$r=1, s=0, \alpha=\beta=\gamma=0, \rho\in${\Bbb R} --arbitrary. 
The results obtained so far are summarized in Table~1.\\
\centerline{\bf Table 1}\samepage
\begin{center}
\begin{tabular}{|c|r|r|r|r|r|r|c|}
\hline
 &$\alpha$&$\beta$&$\gamma$&$\rho$&$r$&$s$&$\mbox{Remarks}$\\ \hline
1&0&0&0&$\varepsilon$&1&0&$\varepsilon \in${\Bbb R} \\ \hline
2&1&$\varepsilon$&$-\varepsilon$&1&0&0&$\varepsilon\geq 0$\\ \hline
3&-1&$\varepsilon$&$-\varepsilon$&-1&0&0&$\varepsilon\geq 0 $\\ \hline
4&-1&$\varepsilon$&$-\varepsilon$&1&0&0&$\varepsilon\geq 0 $\\ \hline 
5&0&0&0&1&0&0& \\ \hline
6&0&1&-1&1&0&0& \\ \hline
7&0&0&0&-1&0&0& \\ \hline
8&0&1&-1&-1&0&0& \\ \hline
9&0&1&-1&0&0&0&\mbox{coboundary} \\ \hline
\end{tabular}
\end{center}
We have checked explictly that all above bialgebra structures are consistent 
and inequivalent. It remains to show that only the last structure is a 
coboundary. Let us remind that a cocommutator \d\  given by

\be
\d(X)=i[1\otimes X+X\otimes 1,\; r],&r\in L\wedge L,&X\in L
\ee

\noindent defines a coboundary Lie bialgebra if and only if $r$ fulfills
the modified classical Yang--Baxter equation

\be
[X\otimes 1\otimes 1+1\otimes X\otimes 1+1\otimes 1\otimes X,\;[[r,r]]]=0,& &\;\;X\in L
\ee

\noindent where $[[r,r]]$\ is the Schouten bracket
$$
[[r,r]]\equiv [r_{12},r_{13}]+[r_{12},r_{23}]+[r_{13},r_{23}];
$$
\noindent here $r_{12}=r^{ij}X_i\otimes X_j\otimes 1$\ etc.

In our case

\be
r&=&aP\wedge H+bP\wedge K+cH\wedge K
\ee

\noindent and

\be
[[r,r]]&=&3ic^2P\wedge K\wedge H
\ee

\noindent and~(22) holds for any values a,b,c; on the other hand classical
Yang--Baxter equation implies $c=0$.\\
Eqs.~(21)~and~(23) give now

\be
\d(P)&=&0\nn\\
\d(H)&=&cH\wedge P\\
\d(K)&=&cK\wedge P\nn
\ee

Under the redefinition \mbox{$H\to H,\; K\to K/c,\;  P\to -P/c$} the above
structure is converted to the one given in the last line of Table~1.

Finally, let us note that the classical r--matrix obeying CYBE defines
here trivial cocomutator.
\section{The Lie--Poisson structures on twodimensional Galilei group}
Let $G$\ be a Lie group, $L$--its Lie algebra,$\{X_i^R\}$--the set
of rightinvariant fields on $G$. Let

\be
\eta(g)&=&\eta^{ij}(g)X_i\otimes X_j
\ee

\noindent be the map $G\to\wedge^2L$. Then

\be
\{\Psi,\Phi\}&\equiv&\eta^{ij}(g)X^R_i\Psi X^R_j\Phi
\ee

\noindent provides $G$\ with a Poisson--Lie group structure
if and only if\\

\noindent $
\mbox{(i)\hspace{0.5cm}}\eta^{il}X^R_i\eta^{jk}+\eta^{kl}X^R_l\eta^{ij}+\eta^{jl}X^R_l\eta{ki}-c^j_{lp}\eta^{il}\eta^{pk}-\\
\mbox{\hspace{1cm}}-c^i_{lp}\eta^{kl}\eta^{pj}-c^k_{lp}\eta^{jl}\eta^{pi}=0 \hfill (28a)\\
$ 

\noindent $
\mbox{(ii)\hspace{0.5cm}}\eta(gh)=\eta(g)+Ad_g\eta(h)\;\hfill (28b);$\\

\noindent $\eta$\ defines the bialgebra structure on $L$\ through
\setcounter{equation}{28}

\be
\d(x)&=&\frac{d\eta(e^{itx})}{dt}|_{t=0}
\ee

Our aim here is to find the Lie--Poisson structures which correspond, 
via~(29), to the bialgebra structures found in previous section. To this
end we use eq.~(28b) to find $\eta$. Then the Poisson bracket is calculated
according to eq.~(26) and Jacobi identities checked. In such a way, following
Ref.~[14] we avoid an explicit use of eqs.~(28a).\\

\ni{Define

\be
\eta(a,v,\tau)&=&\lambda(a,v,\tau)P\wedge H+\mu(a,v,\tau)P\wedge K+\nu(a,v,\tau)H\wedge K
\ee}

\noindent Eq.~(28b) gives\\

\be
&&\eta(a+a^\prime+v\tau^\prime,v+v^\prime,\tau+\tau^\prime)=(\lambda(a,v,\tau)+\lambda(a^\prime,v^\prime,\tau^\prime)+\tau\nu(a^\prime,v^\prime,\tau^\prime))P\wedge H\nn\\
&&\mbox{}+(\mu(a,v,\tau)+\mu(a^\prime,v',\tau')-v\nu(a^\prime,v^\prime,\tau^\prime))P\wedge K\nn\\
&&\mbox{}+(\nu(a,v,\tau)+\nu(a',v',\tau'))H\wedge K
\ee

Consequently, one obtains the following set of equations determining $\lambda,\mu,\nu$

\be
\lambda(a+a'+v\tau', v+v',\tau+\tau')&=&\lambda(a,v,\tau)+\lambda(a',v',\tau')+\tau\nu(a',v',\tau')\nn\\
\mu(a+a'+v\tau',v+v',\tau+\tau')&=&\mu(a,v,\tau)+\mu(a',v',\tau')-v\nu(a',v',\tau')\nn\\
\nu(a+a'+v\tau',v+v',\tau+\tau')&=&\nu(a,v,\tau)+\nu(a',v',\tau')
\ee

\noindent Due to

\be
(a,v,\tau)&=&(0,0,\tau)\star(a,0,0)\star(0,v,0)\nn
\ee

\noindent it is sufficient to find $\eta$\ for one--parameter 
subgroups generated by $H,P$\ and $K$. We get, respectively:\\

\ni$\lambda(0,0,\tau+\tau')=\lambda(0,0,\tau)+\lambda(0,0,\tau')+\tau\nu(0,0,\tau')\\
\mu(0,0,\tau+\tau')=\mu(o,o,\tau)+\mu(0,0,\tau')\hfill (33a)\\
\nu(0,0,\tau+\tau')=\nu(0,0,\tau)+\nu(0,0,\tau')$\\

\ni with\\

\ni$\lambda(0,0,\tau)=\frac{b}{2}\tau^2+c\tau\\
\mu(0,0,\tau)=k\tau\hfill(34a)\\
\nu(0,0,\tau)=b\tau$\\

\ni and

\noindent$\lambda(0,v+v',0)=\lambda(0,v,0)+\lambda(0,v',0)\\
\mu(0,v+v',0)=\mu(0,v,0)+\mu(0,v',0)-v\nu(0,v',0)\hfill (33b)\\
\nu(0,v+v',0)=\nu(0,v,0)+\nu(0,v',0)\\$

\ni with\\

\ni$\lambda(0,v,0)=ev\\
\mu(0,v,0)=fv-\frac{d}{2}v^2\hfill (34b)\\
\nu(0,v,0)=dv\\$

\ni as well as

\ni$\lambda(a+a',0,0)=\lambda(a,0,0)+\lambda(a',0,0)\\
\mu(a+a',0,0)=\mu(a,0,0)+\mu(a',0,0) \hfill (33c)\\
\nu(a+a',0,0)=\nu(a,0,0)+\nu(a',0,0)\\$

\ni with\\

\ni$\lambda(a,0,0)=ga\\
\mu(a,0,0)=ha \hfill (34c)\\
\nu(a,0,0)=ja\\$

Let us now use eqs.~(34) to construct the general form of $\lambda,\mu$\ 
and $\nu$. We write

\setcounter{equation}{34}

\be
(a,v,\tau)&=&((0,0,\tau)\star(a,0,0))\star(0,v,0)
\ee

\noindent and use eqs.~(32) to get

\be
\lambda(a,v,\tau)&=&\frac{b}{2}\tau^2+c\tau+d(v\tau-a)+ev\nn\\
\mu(a,v,\tau)&=&k\tau-ba+fv-\frac{d}{2}v^2\\
\nu(a,v,\tau)&=&b\tau+dv\nn
\ee

Note that the number of free parameters in eqs.~(36) is smaller
than in eqs.~(34). This is due to the fact that once the general forms
of $\lambda,\mu,\nu$\ is calculated from eqs.~(32)~and~(35) they have
to be reinserted back in eqs.~(32) which provide further constraints.
These constraints arise because we have used a specific order of factors
on the right--hand side of eq.~(35) so that the associativity has to be
still imposed.

The general form of $\eta$\ is given by eqs.~(30),(36). However, given 
an automorphism of Galilei group (which results also in some automorphism
of its algebra) we can easily find its action on $\eta$. We call
two $\eta$'s equivalent if they are related by such an automorphism.
Our next aim is to classify nonequivalent $\eta$'s. The simplest way to do
this is to find all maps $\eta$\ giving rise, via eq.~(29), to
cocomutators \d\  classified in previous section. Simple calculation gives

\be
\d(H)&=&-cP\wedge H-kP\wedge K-bH\wedge K\nn\\
\d(P)&=&-dP\wedge H-bP\wedge K\\
\d(K)&=&eP\wedge H+fP\wedge K+dH\wedge K\nn
\ee

\noindent By comparying eqs.~(10)~and~(37) we get

\be
d=-r,&k=\alpha,&f=-\gamma\\
b=s,&c=-\beta,&e=\rho\nn
\ee

Eqs.~(11) impose now further constraints on parameters entering $\eta$.
This is because eqs.~(28a) haven't been used yet.

Table~1 of section~3 can be now used to generate all nonequivalent
Lie--Poisson structures on twodimensional Galilei group. However,
one should keep in mind that the classification in sec.~3 was given for
dimensionless form of Lie algebra while here rather the generators 
with proper dimensions are neccesary (of eqs.~(5),(6)). This
can be cured as described in sec.~3.

Let us first write out the general form of Poisson bracket following
from eqs.~(6),(27) and~(30); it reads\\

\be
\{f,g\}&=&\lambda(\frac{\partial f}{\partial a}\frac{\partial g}{\partial \tau}
-\frac{\partial f}{\partial \tau}\frac{\partial g}{\partial a})+\mu(\frac{\partial f}{\partial v}
\frac{\partial g}{\partial a}-\frac{\partial f}{\partial a}\frac{\partial g}{\partial v})+\nn\\
 & &\mbox{}+\nu(\frac{\partial f}{\partial\tau}(\frac{\partial g}{\partial v}+\tau\frac{\partial g}{\partial a})-
(\frac{\partial f}{\partial v}+\tau\frac{\partial f}{\partial a})\frac{\partial g}{\partial \tau})
\ee

\ni In particular

\be
\{a,v\}=&-\mu&=-k\tau+ba-fv+\frac{d}{2}v^2\nn\\
\{a,\tau\}=&\lambda-\nu\tau&=-\frac{b}{2}\tau^2+c\tau-da+ev\nn\\
\{v,\tau\}=&-\nu&=-b\tau-dv\nn
\ee

Taking care about proper dimension we get finally the set of Lie--Poisson
structures on twodimensional Galilei group corresponding to the bialgebra
structures classified in sec.~3. It is given in table~2\\

\centerline{Table 2}
\begin{center}
\begin{tabular}{|c|c|c|c|c|}
&$\{a,v\}$&$\{a,\tau\}$&$\{v,\tau\}$&Remarks\\ \hline
1&$-\frac{\tau_0 v^2}{2}$&$\tau_0 a+\varepsilon\tau_0^2 v$&$\tau_0 v$&$\varepsilon\in${\Bbb R}\\
2&$-v_0^2\tau-\varepsilon v_0\tau_0v$&$-\varepsilon v_0\tau_0\tau+\tau_0^2v$&0&$\varepsilon \geq 0$\\
3&$v_0^2\tau-\varepsilon v_0\tau_0v$&$-\varepsilon v_0\tau_0\tau-\tau_0^2v$&0&$ \varepsilon \geq 0$\\
4&$v_0^2\tau-\varepsilon v_0\tau_0v$&$-\varepsilon v_0\tau_0\tau+\tau_0^2v$&0&$\varepsilon \geq 0$\\
5&0&$\tau_0^2v$&0&\\
6&$-v_0\tau_0v$&$-v_0\tau_0\tau+\tau_0^2v$&0&\\
7&0&$-\tau_0^2v$&0&\\
8&$-v_0\tau_0v$&$-v_0\tau_0\tau-\tau_0^2v$&0&\\
9&$-v_0\tau_0v$&$-\tau_0v_0\tau$&0&\\ \hline
\end{tabular}
\end{center}

\noindent here $v_0,\tau_0$\ are arbtirary but \underline{fixed} nonzero
constants while $\varepsilon$--dimensionless parameter.
We have checked explictly the Jacobi identities as well as Poisson--Lie
property for all the above cases.

\section{Conclusions}
We have classified above all bialgebra structures on twodimensional 
Galilei algebra and found the corresponding Lie--Poisson structures
on Galilei group. 

It is worthwile to stress the following point. In order to impose
a Lie--Poisson structure on Galilei group one needs two
dimensionful constants. They can attain arbitrary nonzero values,
different choices being related by automorphisms. The only relevant
free parameter is  dimensionless parameter $\varepsilon$;
different values of $\varepsilon$\ correspond to nonequivalent
Lie--Poisson structures.

One should also note that the relatively rich family of 
nonequivalent Lie--Poisson structures considered here contains only
one coboundary; this is in sharp contrast with semisimple case~[15]
as well as the case of fourdimensional Poincare group~[16]~[17].

The more challenging problem is the classification of all Lie--Poisson
structures (resp.~bialgebra structures) on fourdimensional Galilei
group (resp.~algebra). This problem will be considered in forthcoming paper.

\section*{Acknowledgments}
 The author acknowledges prof.~P.~Kosi\'nski a careful reading of the 
manuscript and many helpful suggestion and special thanks to dr~S.~Giller,
dr~C.~Gonera, dr~P.~Ma\'slanka and mgr~A.~Opanowicz for valuable discussion.

\end{document}